\documentclass{article}
\newtheorem{theorem}{Theorem}[section]
\newtheorem{definition}{Definition}[section]

\font\euler=eusm10
\def \M{\mbox{\euler M}}
\def \L{\mbox{\euler L}}

\def \N{\mbox{\euler N}}

\def \intd{{\rm d}}

\def \hbe{{\mbox{\boldmath $\beta$}}}

\def \tx{\tilde{x}}
\def \ty{\tilde{y}}

\def \hhsi{\hat{\sigma}}
\usepackage{amssymb}
\usepackage{amsmath}
\usepackage{bm}
\usepackage{pgfplots}
\usepackage{indentfirst}
\usepackage{multirow}
\usepackage{booktabs}

\begin{document}

\title{Least absolute deviations uncertain regression with imprecise observations}

\author{Zhe Liu$^{*}$ \\
{\em\small Department of Mathematical Sciences, Tsinghua University, Beijing 100084, China}\\
{\em\small z-liu16@mails.tsinghua.edu.cn}\\
{\em\small *Corresponding author}\\
}
\date{} \maketitle

\begin{abstract} 
Traditionally regression analysis answers questions about the relationships among variables based on the assumption that the observation values of variables are precise numbers. It has long been dominated by least squares techniques, mostly due to their elegant theoretical foundation and ease of implementation. However, in many cases, we can only get imprecise observation values and the assumptions upon which the least squares is based may not be valid. So this paper characterizes the imprecise data in terms of uncertain variables and proposes a novel robust approach under the principle of least absolute deviations to estimate the unknown parameters in uncertain regression models. Finally, numerical examples are documented to illustrate our method.

\end{abstract}

\section{Introduction}
As a central part of many research projects, regression analysis is the study of relationships between the response variable and predictor variables by a regression model with a goal to summarize observed data as simply, elegantly, and usefully as possible. Traditionally, regression analysis supposes the observation values of those variables are precise numbers under the framework of probability theory. However in many cases data are imprecisely observed in our daily life and we can not get useful sample data when emergencies such as flood and earthquake occur. Under these situations, many surveys have shown that the probability theory may lead to counterintuitive results and uncertainty theory established by Liu \cite{Liu2007} and refined by Liu \cite{Liu2010} is more suitable for imprecise observations given by experts \cite{Liu2012}. Subsequently many researchers such as Wen et al. \cite{Wen2017}, Lio and Liu \cite{LioL2018, WB}, Nejad and Ghaffari-Hadigheh \cite{NG}, Yao \cite{Yao2018} and Yang and Liu \cite{YangLiu} characterized imprecise observations in terms of uncertain variables in different fields. Especially, uncertain regression analysis estimates the dependence among uncertain variables with imprecisely observed samples. For that matter, Yao and Liu \cite{ure} explored a point estimation for unknown parameters in the model under the principle of least squares. In addition, a prediction interval for the response with new predictor variables in uncertain regression models was suggested by Lio and Liu \cite{WB}. Furthermore, Liu and Jia \cite{LZ} proposed a cross-validation method to evaluate the predictive ability of uncertain regression models.\par
After the model has been defined, the next important task is to estimate the unknown parameters in this model based on the observed data using the chosen estimation method, which refers to as parameter estimation or model fitting. The most widely used and best known method of estimation is called the least squares developed by Legendre and Gauss mainly due to the elegant theoretical foundation and ease of implementation. Under certain assumptions, both the Gauss-Markov theorem and the method of maximum likelihood demonstrate that the least squares is ``best'' with desirable properties. Unfortunately, many authorities think the underlying assumptions upon which the least squares is based may not be valid in practice where this method may result in misleading answers \cite{BS1983}. As a result this situation naturally leads to a requirement that estimation methods should be robust which means results are insensitive to small deviations from the assumptions. That is to say, the method can still maintain good performance when the actual model deviates slightly from theoretical assumptions and a small fraction of the data are altered. Otherwise the good properties of the method under theoretical assumptions have no practical significance.

Among all the approaches to robust regression, the least absolute deviations ($LAD$) regression which was suggested by Boscovich in 1757 and studied by Laplace in 1793 has attracted wide attentions in statistics, engineering, finance, and other fields \cite{Dodge1987}. Although predated than least squares, it was forced into background at first mainly because it has no closed-form solution and must resort to iterative algorithms. Nearly a century later, Edgeworth explored a numerical method to solve the unconstrained $LAD$ problem using the weighted median as a basic operation in each iteration. In order to overcome the cycling problem when dealing with degenerative data \cite{HG1994} in Edgeworth's method, Harris \cite{Harris1950} used linear programming techniques to solve the $LAD$ regression, and Charnes et al. \cite{Charnes1955} minimized the $LAD$ objective function using the simplex method. After that, many simplex-like methods have sprung up where the most representative ones are Barrodale and Roberts \cite{BR1973} and Armstrong et al. \cite{AFK1979}. Other approaches such as the direct decent algorithm suggested by Wesolowsky \cite{Wes1981} and the interior point method proposed by Zhang \cite{Zhang1993} are also efficient. Nowadays $LAD$ estimate can be solved easily by iterative procedures on high speed computers, making it a viable alternative. It is shown that $LAD$ estimate not only has greater power than the least squares estimate for asymmetric error distributions and heavy-tailed, symmetric error distributions but also has greater resistance to the influence of a few outlying values of variables \cite{BD1993}, implying it is actually more efficient in life-like situations where small errors would occur in measurement. In addition to fit regression models, $LAD$ estimate embodies ideas that are important in linear optimization theory and numerical analysis and has been used in other linear situations such as time series and multivariate data analysis, nonlinear regression \cite{WS1994}, classification and regression trees \cite{BFO}, and as a starting estimate for many robust regression methods \cite{Huber}. \par
Note that least squares estimates given in \cite{ure, WB} are vulnerable to outliers. To the best of our knowledge, robust estimate methods with imprecise observations seem not yet be explored. This paper develop and describe a novel robust approach under the principle of least absolute deviations which is resistant to gross deviations of a small number of imprecise observations in the uncertain regression analysis. The rest of this paper is organized as follows. In Section 2, uncertain regression under the principle of least squares is going to be reviewed. A novel least absolute deviations ($LAD$) estimate in uncertain regression models will be introduced in Section 4. After that, Section 5 will give $LAD$ estimates for some specific regression models. In addition, numerical examples in Section 6 are going to show the calculation of the $LAD$ estimate and compare it with the least squares estimate. Finally, Section 7 will conclude the paper with a brief summary. Some basic concepts and properties about uncertainty theory used in this paper will be given in the Appendix.

\section{Uncertain regression with least squares}
In this section, we review the uncertain regression analysis under the principle of least squares. Suppose $(x_{1}, x_{2}, \cdots, x_{p})$ is a vector of predictor variables and $y$ is a response variable. The regression model as a function of the predictor variables $(x_{1},x_{2},\cdots,x_{p})$ is usually formalized as
\begin{equation} \label{basicm2}
  y=g(x_{1},x_{2},\cdots,x_{p}|\hbe)+\epsilon,
\end{equation}
where $\hbe$ is a vector of unknown parameters to be estimated from the observed data, $\epsilon$ is an uncertain error containing information for determining $y$ that is not already captured in $(x_{1},x_{2},\cdots,x_{p})$.\par
Next having observed the data which satisfy the regression model (\ref{basicm2}) we aim to estimate the unknown parameters in this regression model. Traditional regression methods assume that the value of observation data are precise numbers. However in many cases observation values of both $(x_{1},x_{2},\cdots,x_{p})$ and $y$ are imprecise and denoted as
\begin{equation} \label{da1}
(\tx_{i1}, \tx_{i2}, \cdots, \tx_{ip}, \ty_{i}),\quad i=1,2,\cdots,n,
\end{equation}
where $\tx_{i1}, \tx_{i2}, \cdots, \tx_{ip}, \ty_{i}$ are uncertain variables with uncertainty distributions $\Phi_{i1}, \Phi_{i2}, \cdots, \Phi_{ip}, \Psi_{i}$, $i=1,2,\cdots,n$, respectively. Then we have
\begin{equation} \label{residual}
\ty_{i}=g(\tx_{i1}, \tx_{i2}, \cdots, \tx_{ip}|\hbe)+\epsilon_{i}, \quad i=1,2, \cdots, n.
\end{equation}
Under this situation, Yao and Liu \cite{ure} proposed the least squares estimate as the solution of the following minimization problem,
\begin{equation} \label{lse}
\min \limits_{\scriptsize \hbe} \sum \limits_{i=1}^{n} E[(\ty_{i} - g(\tx_{i1},\tx_{i2},\cdots, \tx_{ip}|\hbe))^{2}].
\end{equation}
After that, Lio and Liu \cite{WB} proposed definitions of the $i$-th residuals
\begin{equation} \label{residual1}
\tilde{\epsilon}_{i} = \ty_{i} - g(\tx_{i1}, \tx_{i2}, \cdots, \tx_{ip}|\hbe^{*}),
\end{equation}
$i=1,2, \cdots, n$, respectively, where $\hbe^{*}$ is the least squares estimate of the unknown parameter $\hbe$.
Furthermore if we assume that
\begin{equation*}
  E[\epsilon_{i}]=e, \quad V[\epsilon_{i}]=\sigma^{2},\quad i=1,2, \cdots, n,
\end{equation*}
in Equation (\ref{residual1}), we can use
\begin{equation} \label{expec}
  \hat{e}= \frac{1}{n} \sum \limits_{i=1}^{n}E[\tilde{\epsilon}_{i}]
\end{equation}
to estimate the unknown $e$ and
\begin{equation} \label{varia}
\hhsi ^{2} = \frac{1}{n} \sum \limits_{i=1}^{n}E[(\tilde{\epsilon}_{i} - \hat{e})^{2}]
\end{equation}
to estimate the unknown $\sigma^{2}$.\par
In addition, a regression model is usually constructed for prediction, that is to say, we would like to forecast the value of response variable which has not been observed for future observations of predictor variables based on the given imprecise data (\ref{da1}). Define an uncertain variable $\hat{\epsilon}$ with the expected value $\hat{e}$ and variance $\hat{\sigma}^{2}$, where $\hat{e}$ and $\hat{\sigma}^{2}$ are the estimated expected value and variance of the uncertain error $\epsilon$ in regression model (\ref{basicm2}). Given a vector of new observed data of predictor variables $(\tx_{1}, \tx_{2}, \cdots, \tx_{p})$ independent of $\hat{\epsilon}$, where $\tx_{1}, \tx_{2}, \cdots, \tx_{p}$ are uncertain variables with regular uncertainty distributions $\Phi_{1}, \Phi_{2}, \cdots, \Phi_{p}$, respectively, the forecast uncertain variable $\tilde{\hat{y}}$ of $y$  can be determined as
\begin{equation}\label{fuv}
  \tilde{\hat{y}}=g(\tx_{1}, \tx_{2}, \cdots, \tx_{p}|\hbe^{*})+\hat{\epsilon}.
\end{equation}
Thus as the expected value of forecast uncertain variable $\tilde{\hat{y}}$, the forecast value $\mu$ of $y$ \cite{WB} is
\begin{equation}\label{fv}
  \mu=E[\tilde{\hat{y}}]=E[g(\tx_{1}, \tx_{2}, \cdots, \tx_{p}|\hbe^{*})]+\hat{e}.
\end{equation}
However, as a point prediction the forecast value $\mu$ is too precise to be convincing sometimes while the prediction interval \cite{WB} which has some confidence that our inference must be correct is more suitable to estimate $y$. Taking $\alpha$ as a predetermine level (e.g., 95 \%), we get the $\alpha$ prediction interval of $y$ as
\begin{equation} \label{inter}
[\mu-b, \mu+b]
\end{equation}
in which $b$ is the minimum value such that
\begin{equation*}
\hat{\Psi}(\mu+b)-\hat{\Psi}(\mu-b) \geq \alpha,
\end{equation*}
where $\hat{\Psi}$ is the uncertainty distribution of $\tilde{\hat{y}}$ which can be obtained by the inverse uncertainty distribution $\hat{\Psi}^{-1}$ of $\tilde{\hat{y}}$ and $\mu$ is the forecast value of $y$ given in Equation (\ref{fv}).

\section{Definition of the $LAD$ estimate}
Obviously the least squares estimate is vulnerable to outliers because the square function grows too fast, leading us to wonder whether one can obtain a more robust estimate. Given imprecise observation data (\ref{da1}) which satisfy the regression model (\ref{basicm2}), we shall therefore concentrate our attention to estimates that can be defined by a minimum principle of the form
\begin{equation} \label{rho}
\min \limits_{\scriptsize \hbe} \sum \limits_{i=1}^{n}E\left[ \rho (\ty_{i} - g(\tx_{i1},\tx_{i2},\cdots, \tx_{ip}|\hbe))\right]
\end{equation}
where $\rho$ is some function with the following properties:
\begin{itemize}
  \item $\rho(r) \geq 0$ for all $r$ and has a minimum value $\rho(0)=0$.
  \item $\rho(r)=\rho(-r)$ for all $r$.
  \item $\rho(r)$ increases as $|r|$ increases from $0$.
\end{itemize}
For example, the least squares estimate is a special case by taking $\rho(r)=r^{2}$ in Equation (\ref{rho}). As mentioned earlier, it is usually desirable that the estimate is robust which means outliers in observations do not have unduly large influences on the estimate because errors in observations are inevitable in our real daily life. This produces another property of $\rho(r)$ in Equation (\ref{rho}) as follows:
\begin{itemize}
  \item $\rho(r)$ does not get too large as $r$ increases.
\end{itemize}
As we can see the least square estimate is not robust because $\rho(r)=r^{2}$ does not satisfy the fourth property. Therefore other criteria may be more suitable than least squares to estimate parameters in uncertain regression models when there are some errors in observations. In this section, we discuss a more robust estimate method under the principle of $LAD$ by taking $\rho(r)=|r|$ in Equation (\ref{rho}) to better deal with observations with outliers.
\begin{definition} \label{def1}
  Denote a set of imprecise observation data which satisfy the regression model (\ref{basicm2}) as $(\tx_{i1}, \tx_{i2}$, $\cdots, \tx_{ip}, \ty_{i})$, $i=1,2,\cdots,n$, where $\tx_{i1}, \tx_{i2}$, $\cdots, \tx_{ip}, \ty_{i}$ are independent uncertain variables with regular uncertainty distributions $\Phi_{i1}, \Phi_{i2},$ $\cdots, \Phi_{ip}$,
  $\Psi_{i}$, $i=1,2,\cdots, n$, respectively. So each observation can be written as
   \begin{equation*}
     \ty_{i}=g(\tx_{i1}, \tx_{i2}, \cdots, \tx_{ip}|\hbe)+\epsilon_{i}, \quad i=1,2, \cdots, n.
   \end{equation*}
Then we define the $LAD$ estimate of $\hbe$ as a minimizer of
\begin{equation} \label{basicp}
\min \limits_{\scriptsize \hbe} \sum \limits_{i=1}^{n} E\left| \ty_{i} - g(\tx_{i1},\tx_{i2},\cdots, \tx_{ip}|\hbe) \right|.
\end{equation}
\end{definition}
After getting the $LAD$ estimate $\hat{\hbe}$, the fitted regression model is determined by
\begin{equation} \label{frm}
  y=g(x_{1},x_{2}, \cdots, x_{p}| \hat{\hbe}).
\end{equation}
If we further assume that the function $g(x_{1}, x_{2}, \cdots, x_{p}|\hbe)$ in regression model (\ref{basicm2}) is a strictly monotone function which are satisfied in many practical problems, the minimization problem (\ref{basicp}) in Definition \ref{def1} can be calculated as follows.
\begin{theorem} \label{th41}
Assuming that the function $g(x_{1}, x_{2}, \cdots, x_{p}|\hbe)$ in regression model (\ref{basicm2}) is strictly increasing with respect to $x_{1}, \cdots, x_{m}$ and strictly decreasing with respect to $x_{m+1}, \cdots, x_{p}$, the $LAD$ estimate of $\hbe$ in Definition \ref{def1} can be calculated as
\begin{equation*}
\min \limits_{\hbe} \sum \limits_{i=1}^{n} \int_{0}^{1} \left| \Psi_{i}^{-1}(\alpha) - g(\Phi_{i1}^{-1*}(\alpha), \Phi_{i2}^{-1*}(\alpha), \cdots, \Phi_{ip}^{-1*}(\alpha)|\hbe) \right| \intd \alpha
\end{equation*}
where
\begin{equation*}
\Phi_{ij}^{-1*}(\alpha) =
\left \{ \begin{array}{cc}
\Phi_{ij}^{-1}(1-\alpha), \quad &if \ 1 \leq j \leq m \\
\Phi_{ij}^{-1}(\alpha), \quad &if \ m+1 \leq  j \leq p. \\
\end{array} \right.
\end{equation*}
\end{theorem}

\noindent{\bf Proof:}
Since the function
\begin{equation*}
\ty_{i} -  g(\tx_{i1}, \tx_{i2}, \cdots, \tx_{ip}|\hbe)
\end{equation*}
is strictly increasing with respect to $\ty_{i}$ and strictly decreasing with respect to $\tx_{ij}$ when $1 \leq j \leq m$ or strictly increasing with respect to $\tx_{ij}$ when $m+1 \leq j \leq p$ for each $i$, it follows from Theorem \ref{th1} that the inverse uncertainty distribution is
\begin{equation*}
  F_{i}^{-1}(\alpha)=\Psi_{i}^{-1}(\alpha) - g(\Phi_{i1}^{-1*}(\alpha), \Phi_{i2}^{-1*}(\alpha), \cdots, \Phi_{ip}^{-1*}(\alpha)|\hbe).
\end{equation*}
Then from Equation (\ref{expect}), we obtain
\begin{equation*}
  E\left| \ty_{i} - g(\tx_{i1},\tx_{i2},\cdots, \tx_{ip}|\hbe) \right|=\int_{0}^{1} \left| F_{i}^{-1}(\alpha) \right| \intd \alpha.
\end{equation*}
Thus the minimization problem (\ref{basicp}) can be calculated as
\begin{equation*}
\min \limits_{\hbe} \sum \limits_{i=1}^{n} \int_{0}^{1} \left| \Psi_{i}^{-1}(\alpha) - g(\Phi_{i1}^{-1*}(\alpha), \Phi_{i2}^{-1*}(\alpha), \cdots, \Phi_{ip}^{-1*}(\alpha)|\hbe) \right| \intd \alpha
\end{equation*}
where
\begin{equation*}
\Phi_{ij}^{-1*}(\alpha) =
\left \{ \begin{array}{cc}
\Phi_{ij}^{-1}(1-\alpha), \quad &if \ 1 \leq j \leq m \\
\Phi_{ij}^{-1}(\alpha), \quad &if \ m+1 \leq  j \leq p. \\
\end{array} \right.
\end{equation*}
Then the theorem follows immediately.
\section{$LAD$ estimates for some regression models}
The important instance of regression methodology is linear regression which is the most commonly used in regression analysis. Virtually, many models are generalizations of linear regression models which means they are linear in the unknown parameters after certain transformations, that is, the response variable can be stated in terms of a weighted sum of a set of predictor variables. In fact all other regression methods build upon an understanding of how linear regression works. First we give the $LAD$ estimate in the linear regression model.
\begin{theorem} \label{th2}
Consider the linear regression model in the form
\begin{equation} \label{linear}
\ty_{i} = \beta_{0} + \sum \limits_{j=1}^{p} \beta_{j}\tx_{ij} + \epsilon_{i},\quad i=1, 2, \cdots,n,
\end{equation}
where the imprecise observation data $\tx_{i1}, \tx_{i2},$ $\cdots, \tx_{ip}, \ty_{i}$ are independent uncertain variables with regular uncertainty distributions $\Phi_{i1}, \Phi_{i2},$ $\cdots, \Phi_{ip}$, $\Psi_{i}, i=1,2,\cdots, n$, respectively. Then the $LAD$ estimate of $\hbe=(\beta_{0}, \beta_{1}, \cdots, \beta_{p})$ in Equation (\ref{linear}) solves the following minimization problem:
\begin{equation} \label{mine2}
\min \limits_{\beta_{0},\beta_{1},\cdots,\beta_{p}} \sum \limits_{i=1}^{n}
E \left|\ty_{i} - \beta_{0} - \sum \limits_{j=1}^{p} \beta_{j}\tx_{ij} \right|
\end{equation}
which can be calculated as
\begin{equation*}
\min \limits_{\beta_{0},\beta_{1},\cdots,\beta_{p}} \sum \limits_{i=1}^{n} \int_{0}^{1} \left| \Psi_{i}^{-1}(\alpha) - \beta_{0} - \sum \limits_{j=1}^{p} \beta_{j} \Phi_{ij}^{-1*}(\alpha, \beta_{j}) \right| \intd \alpha
\end{equation*}
where
\begin{equation*}
\Phi_{ij}^{-1*}(\alpha, \beta_{j}) =
\left \{ \begin{array}{cc}
\Phi_{ij}^{-1}(1-\alpha), \quad &if \ \beta_{j} \geq 0 \\
\Phi_{ij}^{-1}(\alpha), \quad &if \ \beta_{j} < 0 \\
\end{array} \right.
\end{equation*}
for $i=1,2,\cdots,n$ and $j=1,2,\cdots,p$.
\end{theorem}

\noindent{\bf Proof:}
According to Definition \ref{def1}, the $LAD$ estimate of $\hbe=(\beta_{0}, \beta_{1}, \cdots, \beta_{p})$ in the linear regression model (\ref{linear}) is actually the optimal solution of the minimization problem,
\begin{equation*}
\min \limits_{\beta_{0},\beta_{1},\cdots,\beta_{p}} \sum \limits_{i=1}^{n}
E \left|\ty_{i} - \beta_{0} - \sum \limits_{j=1}^{p} \beta_{j}\tx_{ij} \right|.
\end{equation*}
Since the function
\begin{equation}
 \beta_{0} + \sum \limits_{j=1}^{p} \beta_{j}\tx_{ij} \nonumber
\end{equation}
is strictly increasing with respect to $\tx_{ij}$ when $\beta_{j} \geq 0$ or strictly decreasing with respect to $\tx_{ij}$ when $\beta_{j}<0$ for each $i$, it follows from Theorem \ref{th41} that the minimization problem (\ref{mine2}) is equivalent to
\begin{equation*}
\min \limits_{\beta_{0},\beta_{1},\cdots,\beta_{p}} \sum \limits_{i=1}^{n} \int_{0}^{1} \left| \Psi_{i}^{-1}(\alpha) - \beta_{0} - \sum \limits_{j=1}^{p} \beta_{j} \Phi_{ij}^{-1*}(\alpha, \beta_{j}) \right| \intd \alpha
\end{equation*}
where
\begin{equation*}
\Phi_{ij}^{-1*}(\alpha, \beta_{j}) =
\left \{ \begin{array}{cc}
\Phi_{ij}^{-1}(1-\alpha), \quad &if \ \beta_{j} \geq 0 \\
\Phi_{ij}^{-1}(\alpha), \quad &if \ \beta_{j} < 0 \\
\end{array} \right.
\end{equation*}
Then the theorem follows immediately. \par
As one of the best-known models of enzyme kinetics, Michaelis-Menten regression model
\begin{equation*}
y = \frac{\beta_{1}x}{\beta_{2}+x} + \epsilon, \quad \beta_{1}>0, \beta_{2}>0
\end{equation*}
describes the rate of enzymatic reactions. The $LAD$ estimate in this model is given as follows.
\begin{theorem} \label{th3}
Consider the Michaelis-Menten regression model in the form
\begin{equation} \label{MM}
\ty_{i} = \frac{\beta_{1}\tx_{i}}{\beta_{2}+\tx_{i}} + \epsilon_{i}, \quad \beta_{1}>0, \beta_{2}>0, \quad i=1, 2, \cdots,n,
\end{equation}
where the imprecise observation data $(\tx_{i}, \ty_{i})$ are independent uncertain variables with regular uncertainty distributions $\Phi_{i}, \Psi_{i}, i=1,2,\cdots, n$, respectively. Then the $LAD$ estimate of $\hbe=(\beta_{1}, \beta_{2})$ in this model
solves the following minimization problem:
\begin{equation} \label{mine3}
\min \limits_{\beta_{1}>0, \beta_{2}>0} \sum \limits_{i=1}^{n}
E \left|\ty_{i} - \frac{\beta_{1}\tx_{i}}{\beta_{2}+\tx_{i}} \right|
\end{equation}
which can be calculated as
\begin{equation*}
\min \limits_{\beta_{1}>0, \beta_{2}>0} \sum \limits_{i=1}^{n} \int_{0}^{1} \left| \Psi_{i}^{-1}(\alpha) - \frac{\beta_{1}\Phi_{i}^{-1}(1-\alpha)}{\beta_{2}+\Phi_{i}^{-1}(1-\alpha)} \right| \intd \alpha.
\end{equation*}
\end{theorem}

\noindent{\bf Proof:}
According to Definition \ref{def1} that the $LAD$ estimate of $\hbe=(\beta_{1}, \beta_{2})$ in the Michaelis-Menten regression model (\ref{MM}) is actually the optimal solution of the minimization problem,
\begin{equation*}
\min \limits_{\beta_{1}>0, \beta_{2}>0} \sum \limits_{i=1}^{n}
E \left|\ty_{i} - \frac{\beta_{1}\tx_{i}}{\beta_{2}+\tx_{i}} \right|.
\end{equation*}
Since the function
\begin{equation*}
\frac{\beta_{1}\tx_{i}}{\beta_{2}+\tx_{i}}
\end{equation*}
is strictly increasing with respect to $\tx_{i}$ for each $i$, it follows from Theorem \ref{th41} that the minimization problem (\ref{mine3}) is equivalent to
\begin{equation*}
\min \limits_{\beta_{1}>0, \beta_{2}>0} \sum \limits_{i=1}^{n} \int_{0}^{1} \left| \Psi_{i}^{-1}(\alpha) - \frac{\beta_{1}\Phi_{i}^{-1}(1-\alpha)}{\beta_{2}+\Phi_{i}^{-1}(1-\alpha)} \right| \intd \alpha.
\end{equation*}
Then the theorem follows immediately. \par
As a sigmoid growth model, the Gompertz regression model
\begin{equation*}
  y=\beta_{1} \exp(-\beta_{2} \exp(-\beta_{3}x)) + \epsilon, \quad \beta_{1}>0, \beta_{2}>0, \beta_{3} >0
\end{equation*}
is especially useful in describing the rapid growth of a certain population of organisms and can account for the eventual horizontal asymptote once the carrying capacity is determined. The $LAD$ estimate in this model is given as follows.
\begin{theorem} \label{th4}
Consider the Gompertz regression model in the form
\begin{equation} \label{G}
\ty_{i}=\beta_{1} \exp(-\beta_{2} \exp(-\beta_{3}\tx_{i})) + \epsilon_{i}, \quad \beta_{1}>0, \beta_{2}>0, \beta_{3} >0, \quad i=1, 2, \cdots,n,
\end{equation}
where the imprecise observation data $(\tx_{i}, \ty_{i})$ are independent uncertain variables with regular uncertainty distributions $\Phi_{i}, \Psi_{i}, i=1,2,\cdots, n$, respectively. Then the $LAD$ estimate of $\hbe=(\beta_{1}, \beta_{2}, \beta_{3})$ in this model solves the following minimization problem:
\begin{equation} \label{mine4}
\min \limits_{\beta_{1}>0, \beta_{2}>0, \beta_{3} >0} \sum \limits_{i=1}^{n}
E \left|\ty_{i} -\beta_{1} \exp(-\beta_{2} \exp(-\beta_{3}\tx_{i})) \right|
\end{equation}
which can be calculated as
\begin{equation*}
\min \limits_{\beta_{1}>0, \beta_{2}>0, \beta_{3} >0} \sum \limits_{i=1}^{n} \int_{0}^{1} \left| \Psi_{i}^{-1}(\alpha) - \beta_{1} \exp(-\beta_{2} \exp(-\beta_{3}\Phi_{i}^{-1}(\alpha))) \right| \intd \alpha.
\end{equation*}
\end{theorem}

\noindent{\bf Proof:}
According to Definition \ref{def1} that the $LAD$ estimate of $\hbe=(\beta_{1}, \beta_{2}, \beta_{3})$ in the Gompertz regression model (\ref{G}) is actually the optimal solution of the minimization problem,
\begin{equation*}
\min \limits_{\beta_{1}>0, \beta_{2}>0, \beta_{3} >0} \sum \limits_{i=1}^{n}
E \left|\ty_{i} -\beta_{1} \exp(-\beta_{2} \exp(-\beta_{3}\tx_{i})) \right|.
\end{equation*}
Since the function
\begin{equation*}
\beta_{1} \exp(-\beta_{2} \exp(-\beta_{3}\tx_{i}))
\end{equation*}
is strictly increasing with respect to $\tx_{i}$ for each $i$, it follows from Theorem \ref{th41} that the minimization problem (\ref{mine4}) is equivalent to
\begin{equation*}
\min \limits_{\beta_{1}>0, \beta_{2}>0, \beta_{3} >0} \sum \limits_{i=1}^{n} \int_{0}^{1} \left| \Psi_{i}^{-1}(\alpha) - \beta_{1} \exp(-\beta_{2} \exp(-\beta_{3}\Phi_{i}^{-1}(\alpha))) \right| \intd \alpha.
\end{equation*}
Then the theorem follows immediately. \par
\section{Simulation studies}
In consideration of the importance of linear regression models, this section employs the model
\begin{equation} \label{nelrm}
  y= \beta_{0}+\beta_{1}x+\epsilon
\end{equation}
to show in detail how to get the $LAD$ estimate of the vector of unknown parameters with imprecise observation data and the robustness of $LAD$ estimate compared with the least squares estimate when there are observation errors.\par

\subsection{Calculation of the $LAD$ estimate}
Suppose $(\tx_{i},\ty_{i})$, $i=1,2,\cdots,15$ in Table $1$ are a set of imprecise observation data which satisfy the linear regression model (\ref{nelrm}), where $\tx_{i},\ty_{i}$ are independent uncertain variables with linear uncertainty distributions, $\Phi_{i},\Psi_{i}$, $i=1,2,\cdots,15$, respectively.
\begin{table}[htbp]
  \centering
  \caption{Imprecise Data where $\L(a,b)$ Represents Linear Uncertain Variable}
  \begin{tabular}{c|c|c}
    \hline
    i & $\tilde{y}_{i}$ & $\tilde{x}_{i}$ \\
    \hline
    1 & \L(2,3) & \L(0,1) \\
    2 & \L(23,24) & \L(7,8) \\
    3 & \L(25,26) & \L(7,8) \\
    4 & \L(7,8) & \L(1,2) \\
    5 & \L(13,14) & \L(3,4) \\
    6 & \L(20,21) & \L(6,7) \\
    7 & \L(31,32) & \L(9,10) \\
    8 & \L(46,47) & \L(15,16) \\
    9 & \L(56,57) & \L(18,19) \\
    10 & \L(74,75) & \L(24,25) \\
    11 & \L(92,93) & \L(30,31) \\
    12 & \L(95,96) & \L(31,32) \\
    13 & \L(38,39) & \L(12,13) \\
    14 & \L(59,60) & \L(19,20) \\
    15 & \L(82,83) & \L(27,28) \\
    \hline
  \end{tabular}
\end{table}
After the data has been collected, the next task is to estimate the unknown parameter $\hbe=(\hbe_{0}, \hbe_{1})$ in the model (\ref{nelrm}) based on the data given in Table 1. In order to get the $LAD$ estimate $\hat{\hbe}=(\hat{\beta_{0}}, \hat{\beta_{1}})$, we solve the minimization problem (\ref{basicp}) according to Definition \ref{def1}, i.e.,
\begin{equation} \label{mp}
\min \limits_{\beta_{0},\beta_{1}} \sum \limits_{i=1}^{15}
E\left| \ty_{i} - (\beta_{0}+\beta_{1}\tx_{i}) \right|.
\end{equation}
Actually, Equation (\ref{mp}) has an equivalent form following from Theorem \ref{th2}, i.e.,
\begin{equation*}
\min \limits_{\beta_{0},\beta_{1}} \sum \limits_{i=1}^{15} \int_{0}^{1} \left| \Psi_{i}^{-1}(\alpha) - \beta_{0} - \beta_{1} \Phi_{i}^{-1*}(\alpha, \beta_{1}) \right| \intd \alpha
\end{equation*}
where
\begin{equation*}
\Phi_{i}^{-1*}(\alpha, \beta_{1}) =
\left \{ \begin{array}{cc}
\Phi_{i}^{-1}(1-\alpha), \quad &if \ \beta_{1} \geq 0 \\
\Phi_{i}^{-1}(\alpha), \quad &if \ \beta_{1} < 0 \\
\end{array} \right.
\end{equation*}
for $i=1,2,\cdots,15$.
Then we obtain the $LAD$ estimate as follows,
\begin{equation*}
  \hat{\hbe}=(\hat{\beta_{0}}, \hat{\beta_{1}})=(2.4016,2.9344).
\end{equation*}
As a result the fitted linear regression model is
\begin{equation*}
  y=2.4016+2.9344x.
\end{equation*}
It follows from Theorems \ref{th1} and Equation \ref{expect} that Equation (\ref{expec}) has an equivalent form,
\begin{equation*}
  \hat{e}=\frac{1}{15} \sum_{i=1}^{15} \int_{0}^{1} \left( \Psi_{i}^{-1}(\alpha)-2.4016-2.9344\Phi_{i}^{-1}(1-\alpha) \right) \intd \alpha,
\end{equation*}
and Equation (\ref{varia}) has an equivalent form,
\begin{equation*}
  \hat{\sigma}^{2}=\frac{1}{15} \sum_{i=1}^{15} \int_{0}^{1} \left( \Psi_{i}^{-1}(\alpha)-2.4016-2.9344\Phi_{i}^{-1}(1-\alpha)-\hat{e} \right)^{2} \intd \alpha.
\end{equation*}
As a result the estimate expected value $\hat{e}$ and variance $\hat{\sigma}^{2}$ of $\epsilon$ are
\begin{equation*}
  \hat{e}=-0.0548, \quad \hat{\sigma}^{2}=1.3689,
\end{equation*}
respectively.
Suppose we have a new observation of the uncertain predictor variable with $\tx \sim \L(5,6)$ which is independent of $\hat{\epsilon}$, where $\hat{\epsilon}$ is an uncertain variable with the expected value $\hat{e}$ and variance $\hat{\sigma}^{2}$. According to Equation (\ref{fuv}) the forecast uncertain variable of the response variable $y$ is
\begin{equation*}
  \tilde{\hat{y}}=2.4016+2.9344\tx+\hat{\epsilon}.
\end{equation*}
Then the point forecast $\mu$ of $y$ is $18.485$ calculated by Equation (\ref{fv}), i.e.,
\begin{equation*}
  \mu=E[\tilde{\hat{y}}]=2.4016+2.9344E[\tx]+\hat{e}.
\end{equation*} \par
Assuming that $\hat{\epsilon}$ is a normal uncertain variable $\N(\hat{e}, \hhsi)$, we want to get the prediction interval of $y$ with confidence level $\alpha=90 \%$. The uncertainty distribution $\hat{\Psi}$ of $\tilde{\hat{y}}$ can be obtained by the inverse uncertainty distribution $\hat{\Psi}^{-1}$, i.e.,
 \begin{equation}
\hat{\Psi}^{-1}(\alpha) = 2.4016 + 2.9344 (5(1-\alpha) + 6\alpha)+ \Upsilon^{-1}(\alpha)
\end{equation}
where $\Upsilon^{-1}(\alpha)$ is the inverse uncertainty distribution of $\hat{\epsilon}$ with
\begin{equation*}
\Upsilon^{-1}(\alpha)=\hat{e}+ \frac{\hhsi \sqrt{3}}{\pi} \ln \frac{\alpha}{1-\alpha}.
\end{equation*}
The minimum value of $b$ such that
\begin{equation*}
\hat{\Psi}(\mu+b)-\hat{\Psi}(\mu-b) \geq 90 \%
\end{equation*}
is $3.2198$.
It follows from Equation (\ref{inter}) that the $90\%$ prediction interval of the response variable $y$ is
\begin{equation*}
 [15.2652, 21.7948].
\end{equation*}
\par
\subsection{Comparison between $LAD$ and least squares}
Then we conduct another numerical example to show the robustness of $LAD$ estimate compared with least squares estimate when there are observation errors in imprecisely observed data of uncertain variables. Such errors are often encountered within a firm in which highly disaggregated data are used.\par
In the linear regression model (\ref{nelrm}) the parameter $\hbe=(\beta_{0}, \beta_{1})$ is defined as $\hbe=(10,2)$ for simulation purpose, i.e.,
\begin{equation*}
  y=10+2x+\epsilon,
\end{equation*}
and three sets of the imprecisely observed data $(\tx_{ji},\ty_{i})$, $j=1, 2, 3$, $i=1, 2, \cdots ,10$, respectively, are given in Table 2.
 \begin{table}[htbp]
  \centering
  \caption{Imprecise Data where $\L(a,b)$ Represents Linear Uncertain Variable}
  \begin{tabular}{c|c|c|c|c}
    \hline
      &  Models   & $1$ & $2$ & $3$ \\
      \hline
    i & $\tilde{y}_{i}$ & $\tilde{x}_{1i}$ & $\tilde{x}_{2i}$ & $\tilde{x}_{3i}$ \\
    \hline
    1 & \L(10,12) & \L(0,1) & \L(0,1) & \L(0,1)\\
    2 & \L(14,16) & \L(2,3) & \L(5,6)& \L(2,3)\\
    3 & \L(18,20) & \L(4,5) & \L(4,5)& \L(20,21)\\
    4 & \L(22,24) & \L(6,7) & \L(6,7)& \L(6,7)\\
    5 & \L(26,28) & \L(8,9) & \L(8,9)& \L(8,9)\\
    6 & \L(30,32) & \L(10,11) & \L(10,11)& \L(10,11)\\
    7 & \L(34,36) & \L(12,13) & \L(12,13)& \L(12,13)\\
    8 & \L(38,40) & \L(14,15)& \L(14,15)&  \L(14,15)\\
    9 & \L(42,44) & \L(16,17) & \L(6,7) & \L(16,17)\\
    10 & \L(46,48) & \L(18,19)& \L(18,19)&  \L(8,9)\\
    \hline
  \end{tabular}
\end{table}
As we can see, model $1$ with $j=1$ contains no observation errors in $(\tx_{1i}, \ty_{i})$, $i=1, 2, \cdots, 10$, model $2$ with $j=2$ contains observation errors in $(\tx_{2i}, \ty_{i})$ with $2$nd and $9$th observed values, and model $3$ with $j=3$ contains observation errors in $(\tx_{3i}, \ty_{i})$ with $3$rd and $10$th observed values. The corresponding $LAD$ estimates $\hat{\hbe}^{j}$ and least squares estimates $\hbe^{*j}$ for $\hbe$ in model $j$, $j=1, 2, 3$, respectively, are both given in Table 3.

\begin{table}[htbp]
  \centering
  \caption{Estimates of $\hbe$ where $\hat{\hbe}^{j}$ are the $LAD$ estimates and $\hbe^{*j}$ are the least squares estimates, $j=1, 2, 3$, respectively}
  \begin{tabular}{c|c|c}
    \hline
    j & $(\hat{\beta}^{j}_{0}, \hat{\beta}^{j}_{1})$ & $(\beta^{j*}_{0}, \beta^{j*}_{1})$ \\
    \hline
    1 &  (10, 2) &  (10.0479, 1.995)\\
    \hline
    2 & (10, 2) &  (12.3695, 1.8898) \\
    \hline
    3 & (10, 2) & (19.5837, 0.9323) \\
    \hline
  \end{tabular}
\end{table}
Next we analyze the result in Table $3$. As would be expected, in model $1$ which has no observation errors in the observed data, both the least squares estimate $\hbe^{*1}=(10.0479, 1.995)$ and $LAD$ estimate $\hat{\hbe}^{1}=(10,2)$ are very close in value to the true value $\hbe=(10,2)$. The estimation accuracy is almost the same between two methods under this situation. In model $2$ with observation errors in the given imprecisely observed data the $LAD$ estimate $\hat{\hbe}^{2}=(10, 2)$ still equals to the true value $\hbe$ while the least squares estimate $\hbe^{*2}=(12.3695, 1.8898)$ is more far away from the true value $\hbe$ compared with the $LAD$ estimate. In model $3$ the superiority of $LAD$ estimate $\hat{\hbe}^{3}=(10,2)$ is more obvious compared with least squares estimate $\hbe^{*3}=(19.5837, 0.9323)$. \par
Then we delete the outliers in models $2$ and $3$. That is to say, we delete the 2nd data $(\tx_{22}, \ty_{2})$ and the 9th data $(\tx_{29}, \ty_{9})$ in model 2, and delete the 3rd data $(\tx_{33}, \ty_{3})$ and the 9th data $(\tx_{39}, \ty_{9})$ in model 3. With the remaining observations, the $LAD$ estimates $\hat{\hbe}_{-}^{j}$ and least squares estimates $\hbe_{-}^{*j}$ for $\hbe$ in model $j$, $j=2,3$, respectively, are shown in Table $4$.
\begin{table}[htbp]
  \centering
  \caption{Estimates after deleting the outliers}
  \begin{tabular}{c|c|c}
    \hline
    j & $(\hat{\beta}^{j}_{-0}, \hat{\beta}^{j}_{-1})$ & $(\beta^{*j}_{-0}, \beta^{*j}_{-1})$ \\
    \hline
    2 & (10, 2) &  (10.1089, 1.9885) \\
    \hline
    3 & (10.16, 1.9821) & (10.1078, 1.9880) \\
    \hline
  \end{tabular}
\end{table}\par
We can see in Table $4$, the least squares method gives more reasonable results by removing the outliers from the fitting equations.  However, outliers are difficult to distinguish first in many cases. On the other hand, the $LAD$ method achieve the same result without the procedure to exclude outliers by providing residuals contaminated less by the effects of the anomalous observations. \par
Results in above simulations demonstrate that $LAD$ estimate is more robust than least squares with outliers in the given imprecisely observed data, which implies that the estimate under the principle of $LAD$ is more suitable in life-like situations where small errors in observations are inevitable. In fact when abnormal data exist, the least squares reduce the abnormal degree of the abnormal data at the expense of the fitting degree of normal data. Obviously it is harmful because it conceals the truth which means some robust regression method such as $LAD$ method is more appropriate than the least squares to deal with observations with outliers.

\section{Conclusion}
Actually we are usually in the situation that the observation data are imprecise which can not be denoted as precise numbers. What's more, some mistakes are inevitable when collecting data, resulting in outliers in the observation data. Under this situation, classical regression analysis and least squares may lead to counterintuitive results. So this paper introduced $LAD$ estimate in uncertain regression analysis to handle the imprecise data with outliers reasonably. Furthermore, numerical examples were documented to show the calculation for unknown parameters under the principle of $LAD$ and the robustness of the $LAD$ estimate compared with least squares estimate, which showed that $LAD$ estimate is actually more efficient when there are errors in observations.

\section*{Acknowledgments}

\section*{Appendix}\label{appendix}

The fundamental assumption of probability theory is that we can obtain a probability distribution which is close enough to the frequency of the indeterminate quantity. Unfortunately, in practice this assumption may not be valid in all. Under this situation in order to better deal with indeterminacy, Professor Baoding Liu \cite{Liu2007} established the uncertainty theory in 2007. First we review some fundamental concepts, properties and theorems in uncertainty theory. \par
Assume that $\Gamma$ is a nonempty set, and $\L$ is a $\sigma$-algebra over $\Gamma$. Each element $\Lambda$ in $\L$ is called an event. Uncertain measure $\M$ defined by Liu \cite{Liu2007} which indicates the belief degree that an uncertain event may happen satisfies the following three axioms:

\noindent Axiom 1. (Normality Axiom) $\M \{ \Gamma \}=1$ for the universal set $\Gamma$.

\noindent Axiom 2. (Duality Axiom) $\M \{\Lambda\}+\M \{\Lambda^{c}\}=1$
for any event $\Lambda$.

\noindent Axiom 3. (Subadditivity Axiom) For every countable sequence of events $\Lambda_1, \Lambda_2, \!\cdots,$ we have
\[\displaystyle\M \left\{\bigcup_{i=1}^{\infty}\Lambda_i\right\}\le\sum_{i=1}^{\infty}\M \{\Lambda_i\}.\]
Furthermore, Liu \cite{Liu2009} defined the product uncertain measure on the product $\sigma$-algebra $\L$ producing the fourth axiom of uncertainty theory.\\
\noindent Axiom 4. (Product Axiom) {\rm (Liu \cite{Liu2009})} Let $(\Gamma_k,\L_k,\M_k)$ be uncertainty spaces for $k=1, 2, \cdots .$ The product uncertain measure $\M$ is an uncertain measure satisfying
\[\displaystyle\M \left\{\prod_{k=1}^\infty\Lambda_k\right\}=\bigwedge_{k=1}^\infty\M_k\{\Lambda_k\},\]
where $\Lambda_k$ are arbitrarily chosen events from $\L_k$ for $k=1, 2, \cdots$, respectively.
\begin{definition} {\rm (Liu \cite{Liu2007})}
Let $\Gamma$ be a nonempty set, let $\L$ be a $\sigma$-algebra over $\Gamma$, and let $\M$ be an uncertain measure. Then the triplet $(\Gamma,\L,\M)$ is called an uncertainty space.
\end{definition}\par
In addition, Liu \cite{Liu2007} proposed the concept of uncertain variable to represent quantities with uncertainty and the concept of uncertainty distribution to describe uncertain variables. An uncertain variable $\xi$ is a measurable function from the uncertainty space $(\Gamma, \L, \M)$ to the set of real numbers such that for any Borel set B of real numbers, the set
\begin{equation*}
\{ \xi \in B \} = \{ \gamma \in \Gamma \ | \ \xi (\gamma) \in B \}
\end{equation*}
is an event. The uncertain variables $\xi_{1}, \xi_{2}, \cdots, \xi_{n}$ are said to be independent \cite{Liu2009} if
\[\displaystyle\M \left \{ \bigcap_{i=1}^{n} (\xi_{i} \in B_{i}) \right \} = \bigwedge_{i=1}^{n} \M \left \{ \xi_{i} \in B_{i} \right \} \]
for any Borel sets $B_{1}$, $B_{2}$, $\cdots$, $B_{n}$ of real numbers.
 The uncertainty distribution $\Phi$ of an uncertain variable $\xi$ is defined by
\begin{equation}
\Phi (x) =  \M \{ \xi \leq x \} \nonumber
\end{equation}
for any real number $x$. An uncertainty distribution $\Phi (x)$ is said to be regular if it is a continuous and strictly increasing function with respect to $x$ at which $0<\Phi(x)<1$, and
\[\displaystyle \lim_{x \to -\infty} \Phi(x) = 0, \qquad \lim_{x \to \infty} \Phi(x) = 1. \]
Assuming that $\xi$ is an uncertain variable with regular uncertainty distribution $\Phi(x)$, the inverse function $\Phi^{-1}(\alpha)$ is called the inverse uncertainty distribution \cite{Liu2007} of $\xi$.
Generally speaking, the inverse uncertainty distribution for a strictly monotone function of independent uncertain variables can be obtained as follows.
\begin{theorem}{\rm (Liu \cite{Liu2010})} \label{th1}
Let $\xi_{1}, \xi_{2}, \cdots, \xi_{n}$ be independent uncertain variables with regular uncertainty distributions $\Phi_{1}, \Phi_{2}, \cdots, \Phi_{n}$, respectively. If $f$ is strictly increasing with respect to $\xi_{1}, \xi_{2}, \cdots, \xi_{m}$ and strictly decreasing with respect to $\xi_{m+1}, \xi_{m+2}, \cdots, \xi_{n}$, then $\xi = f(\xi_{1}, \xi_{2}, \cdots, \xi_{n})$ is an uncertain variable with an inverse uncertainty distribution
\[\displaystyle \Psi^{-1}(\alpha) = f(\Phi^{-1}_{1}(\alpha), \cdots, \Phi^{-1}_{m}(\alpha), \Phi^{-1}_{m+1}(1-\alpha), \cdots, \Phi^{-1}_{n}(1-\alpha)).\]
\end{theorem}
\begin{definition} {\rm (Liu \cite{Liu2007})}
  Let $\xi$ be an uncertain variable. Then the expected value of $\xi$ is defined as
  \begin{equation}\label{expect1}
   E[\xi] = \int_{0}^{+\infty} \M \{\xi \geq x \}\mathrm{d} x - \int_{-\infty}^{0} \M \{\xi \leq x \}\mathrm{d} x
  \end{equation}
   provided that at least one of the two integrals is finite.
\end{definition}
\begin{theorem} {\rm (Liu \cite{Liu2010})}
  Let $\xi$ be an uncertain variable with regular uncertainty distribution $\Phi$. Then we have
  \begin{equation} \label{expect}
    E\left[ \xi \right]=\int_{0}^{1} \Phi^{-1}(\alpha) \intd \alpha.
  \end{equation}
  \begin{equation}\label{absolute}
  E \left| \xi \right|  = \int_{0}^{1} \left| \Phi^{-1}(\alpha) \right| \intd \alpha.
\end{equation}
\begin{equation} \label{2moment}
    E\left[ \xi^{2} \right]=\int_{0}^{1} (\Phi^{-1}(\alpha))^{2} \intd \alpha.
  \end{equation}
\end{theorem}

\end{document}